\documentclass[11pt,showlabels]{article1}
\usepackage{amssymb}
\usepackage{epsfig}
\newtheorem{theorem}{Theorem}
\newtheorem{assumption}{Assumption}
\newtheorem{lemma}{Lemma}
\newtheorem{remark}{Remark}
\newtheorem{definition}{Definition}

\title{\bf Optimal Stabilization Control for Discrete-time Markov Jump Linear System with Control Input Delay\thanks{This work is supported by the National Natural Science
Foundation of China (Nos. 61473134, 61573220, 61120106011, 61573221) and the Postdoctoral Science Foundation of China (No.
2017M622231). $^{*}$Corresponding author: Huanshui Zhang. Email: hszhang@sdu.edu.cn}
}
\author{Chunyan Han$^a$,\ Hongdan Li$^b$,\ Huanshui Zhang$^{b,*}$
\ \\
\\
\ \ \ $^a$ School of Electrical Engineering, University of Jinan, \\Jinan Shandong 250022, China
\\
$^b$ School of Control Science and Engineering, Shandong
University, \\Jinan Shandong 250061, China}

\begin{document}
\baselineskip 16pt
\date{}
  \maketitle
\begin{abstract}
This paper will investigate the infinite horizon optimal control and stabilization problems for the
Markov jump linear system (MJLS) subject to control input delay. Different from previous
works, for the first time, the necessary and sufficient stabilization
conditions are explored under explicit expressions, and the optimal
controller for infinite horizon is designed with a coupled
algebraic Riccati equation. By
introducing a new type of Lyapunov equation, we show that
under the exact observability assumption, the
MJLS with control input delay is stabilizable in the mean square sense with
the optimal controller if and only if a coupled algebraic Riccati equation has a unique
positive definite solution. The presented
results are parallel to the optimal control and stabilization for standard system with input delay.

\bigskip

\noindent \textbf{Keywords:} Optimal control, stabilization, algebraic Riccati equation, Markov jump linear system, input delay.
\end{abstract}

\pagestyle{plain} \setcounter{page}{1}
\section{Introduction}

Discrete-time MJLSs represent an important class of stochastic systems because they can be used to model random abrupt changes in structure. Dynamic systems with delays \cite{NBW98} or more general networked control applications \cite{HNX07}, where communication networks are used to interconnect remote sensors, actuators and processors, have been shown to be amendable to MJLS modeling. Motivated by a wide spectrum of applications, there has been active research in the stabilization control problems for state delayed MJLSs\cite{XL06}-\cite{SWSSK19}, while no optimality of the controller was considered simultaneously. Different from the previous stabilization results, we mainly study  the stabilization properties of the infinite horizon optimal controller in this paper.

The optimal control and stabilization problems for infinite horizon MJLSs have attracted the attention of many researchers. For example,  necessary and sufficient existence conditions for the infinite horizon optimal controller were developed in \cite{CWC86}, and sufficient stabilization conditions were explored which containing infinite sums. In \cite{JC88}, the definitions of the controllability and observability of discrete-time MJLS were developed. The existence of optimal steady-state controllers was guaranteed by the absolute controllability and the stability of the controlled system was guaranteed by the absolute observability. In \cite{AFJ95}, a necessary and sufficient condition was presented for the existence of a positive-semidefinite solution of the coupled algebraic Riccati-like equation (CARE) occurring in the infinite horizon optimal control problems. In \cite{CFM05}, the concepts of mean square stabilizability and detectability were proposed, and a necessary and sufficient condition for the existence of the stabilizable solution to the infinite horizon optimal control problem was proposed under these concepts. In \cite{CV01}, a new detectability concept (weak detectability) for discrete-time MJLS was presented, and the new concept supplied a sufficient condition for the mean square stability of the infinite-horizon optimal controlled system. It has been shown in \cite{CV01} that mean square detectability developed in \cite{CFM05} ensured weak detectability. Latter, Costa and do Val \cite{CV02} summarized the available results and gave a proposition on the mean square stabilizable of the system. Under the assumption that the system was weak detectable, the system was mean square stabilizable if and only if there existed a positive semi-definite solution to the CARE. More recently, the optimal stationary output feedback control for MJLS waere studied in \cite{CF16}, \cite{DC16}. It can be seen that no time delays are considered in the aforementioned results.

The optimal control and stabilization of
linear system with input delay has received renewed interest in recent years. Various approaches have been developed. The prediction method is one of the famous methods to deal with input delays, which can be traced to the Smith predictor \cite{Smith59}. To overcome the limitation of the original Smith predictor that is just suitable for the open-loop stable system, the finite
spectrum assignment and reduction method are developed in \cite{MO79} and \cite{Art82}, respectively. It can be found that traditional prediction method contains the computation of distributed terms, which is computationally challenging and may be source of instability. To overcome
this problem, the truncated predictor feedback \cite{ZLD12} and closed-loop predictor approaches \cite{ZLM17}, \cite{CG17}, have been developed. Obviously, prediction method is an efficient tool for dealing with the optimal control and stabilization problems for the linear systems with input delay. However, it can't be applied to the stochastic systems subject to input delay directly since the separation principle is not satisfied in the stochastic system. In \cite{ZLXF15} and \cite{LZ16}, the stabilization problems for the discrete-time systems with multiplicative noise and input delays are considered. It has been shown that the system under consideration can be stabilized if and only if the algebraic
Riccati-type equations have a unique solution such that a specific matrix is positive definite.
It is worth mentioning that few results
about infinite horizon optimal control and stabilization problem for MJLS with input delay have been published, which motivates us to undertake the present study. The fundamental questions we will answered in this paper are that: 1) To explore the necessary and sufficient condition for the existence of the infinite horizon optimal controller for the MJLS with input delay; 2) To find the mean square stablization condition for the optimal controlled MJLS.

 In this paper, we aim to provide a thorough solution to the
problems of optimal control and stabilization for infinite
horizon discrete MJLS subject to input delay. As the preliminaries,
the results of finite horizon optimal control for MJLS with control input delay are presented,
and the necessary and sufficient solvability condition of
finite horizon case is given in an explicit expression. By doing
the convergence analysis on the coupled Riccati equation for
the finite horizon case, the infinite horizon optimal controller
and the stabilization conditions (necessary and sufficient) are
derived. In addition, the Lyapunov function for stabilization is
expressed with the optimal cost function. The stabilization
result is obtained under the assumption of exact
observability, under which it is shown that the optimal controlled MJLS
 is mean square stabilizable if and only if the coupled
ARE has a unique positive definite solution.

The remainder of this paper is organized as follows. Section
2 provides the problem formulation and the preliminary results
of finite horizon optimal control for MJLS with input delay. In Section 3, main
results of the infinite horizon optimal control and stabilization
problems are presented. Numerical examples are given in
Section 4 to illustrate main results of this paper. Some
concluding remarks are given in Section 5. Finally, relevant
proofs are detailed in Appendices.

Notations: Throughout this paper, ${R}^n$ denotes the $n$-dimensional Euclidean space,
$R^{m\times n}$ denotes the norm bounded linear space of all
$m\times n$ matrices.
For $L\in R^{n\times n}$, $L'$ stands for the transpose of $L$. As usual, $L\geq 0 (L>0)$
will mean that the symmetric matrix $L\in R^{n\times n}$ is positive semi-definite
(positive definite), respectively.
 $\mbox{E}(.)$ denotes the mathematical expectation operator, $\mbox{P}(.)$ means the occurrence probability of an event.
$\{\Omega,{\cal{G}}, {\cal{G}}_k, P\}$ represents a stochastic basis, with ${\cal{G}}_k$ the $\sigma$-field generated by the random variables $\{x(s),\theta(s); s=0,1,\cdots,k\}$.
We will compactly write the sum $\sum_{l_{k-d+1}=1}^L\lambda_{l_{k-d},l_{k-d+1}}~\cdots ~\sum_{l_{k}=1}^L\lambda_{l_{k-1},l_k}$ as $\Lambda_{l_{k-d},l_{k}}$ and denote $\sum_{l_{d-1}=1}^L\lambda_{l_{d},l_{d-1}}~\cdots ~\sum_{l_{0}=1}^L\lambda_{l_{1},l_0}$ as $\Lambda_{l_{d},l_{0}}$.


\section{Problem Formulation and Preliminaries}

\subsection{Problem Formulation}

\setcounter{equation}{0}
We consider in this paper the infinite horizon optimal control and stabilization problems for the MJLS with input delay.
On the stochastic basis $(\Omega,{\cal{G}}, {\cal{G}}_k,\mbox{P})$, consider the following time-invariant MJLS
\begin{eqnarray}
x(k+1)&=&A_{\theta(k)}x(k)+B_{\theta(k)}u(k-d),\label{fff1}
\end{eqnarray}
where $x(k)\in {\mbox{R}}^n$ is the state, $u(k)\in {\mbox{R}}^m$  is the control input with delay $d>0$. $\theta(k)$ denotes the system mode which
is a discrete-time homogeneous Markov chain. It is assumed that $\theta(k)$ takes values in a finite state space $\Theta\triangleq\{1,2,\cdots,L\}$ with transition probability matrix $\Pi\triangleq (\lambda_{ij})$, where $\lambda_{ij}\triangleq P(\theta(k+1)=j|\theta(k)=i)\geq 0$ for $i,j\in \Theta$ and $\sum_{j=1}^L\lambda_{ij}=1$ for every $i\in \Theta$. The matrices $A_{\theta(k)}$ and $B_{\theta(k)}$ are selected at each time step $k$ from time invariant sets $\{A_1,\cdots,A_L\}$ and $\{B_1,\cdots,B_L\}$ according to the value of the mode. We assume that $\theta(k)$ is independent of $x_0$ and $u(i),i=-d,\cdots,-1$, and the initial values $x_0,u(i),i=-d,\cdots,-1$ are known.

The associated cost function is defined as
\begin{eqnarray}
J\triangleq\mbox{E}\{\sum_{k=0}^\infty x(k)'Qx(k)+\sum_{k=d}^{\infty}u(k-d)'Ru(k-d)\},\label{ss1}
\end{eqnarray}
where $Q$ and $R$ are deterministic symmetric weighting matrices with appropriate dimensions.

Now, we present some definitions.
\begin{definition}
The discrete-time MJLS (\ref{fff1}) is said to be mean square stabilizable if for any initial conditions, there exist a ${\cal{G}}_k$-measurable constant gain controller
\begin{eqnarray}
u(k)=-K_{\theta(k)}^0x(k)-\sum_{i=1}^{d-1}K_{\theta(k)}^iu(k+i-d)\nonumber
\end{eqnarray}
such that $\lim _{k\rightarrow \infty}\mbox{E}(x(k)'x(k))=0$.
\end{definition}
\begin{definition}
The following MJLS
\begin{eqnarray}
x(k+1)=A_{\theta(k)}x(k),y(k)=C_{\theta(k)}x(k),\label{ss2}
\end{eqnarray}
\end{definition}
is called exactly observable, if for any $N\geq n$,
\begin{eqnarray}
y(k)\equiv 0, a.s., \forall 0\leq k\leq N\Rightarrow x_0=0.\nonumber
\end{eqnarray}
For simplicity, we rewrite system (\ref{ss2}) as $(\bar{A},\bar{C})$, where $\bar{A}=(A_1,\cdots,A_{L}), \bar{C}=(C_1,\cdots,C_L)$ with $C_i=C, (i=1,\cdots,L)$, and $Q=C'C$.

Throughout this paper, the following assumptions are required to be satisfied.
\begin{assumption}
$R>0$ and $Q\geq 0$.
\end{assumption}

\begin{assumption}
$(\bar{A},\bar{C})$ is exactly observable.
\end{assumption}

Then the infinite-horizon optimal control and stabilization for system (\ref{fff1}) can be stated as:

\textbf{Problem 1}: Find a ${\cal{G}}_k$-measurable controller $u(k)=-K_{\theta(k)}^0x(k)-\sum_{i=1}^{d-1}K_{\theta(k)}^iu(k+i-d)$, which minimizes the cost (\ref{ss1}) and stabilizes system (\ref{fff1}).

\subsection{Preliminaries}
In this part, we recall the finite horizon optimal control problem for system (\ref{fff1}), which serves as the preparation for the infinite horizon case.

The finite horizon cost function
associated to system (\ref{fff1}) with admissible control law $u=(u(0),\cdots,u(N-d))$ is given by
\begin{eqnarray}
J_N&=&\mbox{E}[\sum_{k=0}^Nx(k)'Qx(k)+\sum_{k=d}^Nu(k-d)'Ru(k-d)\nonumber\\
&&+x(N+1)'P_{\theta(N+1)}x(N+1)],\label{f2}
\end{eqnarray}
where $N>d$ is an integer, $Q, R$ and $P_{\theta(N+1)}$ are deterministic symmetric matrices with compatible dimension and $P_{\theta(N+1)}\geq 0$.

The finite-horizon optimal control for system (\ref{fff1}) can be stated as: Find a ${\cal{G}}_k$-measurable $u(k)$ such that (\ref{f2}) is minimized subject to (\ref{fff1}).

For the convenience of discussions latter, we denote a realization (sample path) of $\{\theta(k-d),\cdots,\theta(k)\}$ by $\{l_{k-d},\cdots,l_k\}$, and define a backward coupled difference equation
\begin{eqnarray}
W_{l_{k-d}}(k-d)&=&\Lambda_{l_{k-d},l_k}[B_{l_k}'(P_{l_k}(k)-P_{l_k}^0(k))B_{l_k}+R] \nonumber\\
&&-\sum_{s=1}^{d-1}\{\Lambda_{l_{k-d},l_{k-s}}[(T_{l_{k-s}}^s(k-s))'W_{l_{k-s}}(k-s)^{-1}
T_{l_{k-s}}^s(k-s)]\},\label{a1}\\
T_{l_{k-d}}^0(k-d)&=&\Lambda_{l_{k-d},l_k}[B_{l_k}'(P_{l_k}(k)-P_{l_k}^0(k))F_{l_k,l_{k-d+1}}]\nonumber\\
&&-\sum_{s=1}^{d-1}\{\Lambda_{l_{k-d},l_{k-s}}[(T_{l_{k-s}}^s(k-s))'W_{l_{k-s}}(k-s)^{-1}T_{l_{k-s}}^0(k-s)\nonumber\\
&&\times F_{l_{k-s},l_{k-d+1}}]\},\label{a2}\\
T_{l_{k-d}}^1(k-d)&=&\Lambda_{l_{k-d},l_k}[B_{l_k}'(P_{l_k}(k)-P_{l_k}^0(k))F_{l_{k},l_{k-d+2}}B_{l_{k-d+1}}]\nonumber\\
&&-\sum_{s=1}^{d-2}\{\Lambda_{l_{k-d},l_{k-s}}[(T_{l_{k-s}}^s(k-s))'W_{l_{k-s}}(k-s)^{-1}T_{l_{k-s}}^0(k-s)\nonumber\\
&&\times F_{l_{k-s},l_{k-d+2}}B_{l_{k-d+1}}]\}-\Lambda_{l_{k-d},l_{k-d+1}}\nonumber\\
&&\times [(T_{l_{k-d+1}}^{d-1}(k-d+1))'W_{l_{k-d+1}}(k-d+1)^{-1}T_{l_{k-d+1}}^0(k-d+1)B_{l_{k-d+1}}],\label{a3}\\
T_{l_{k-d}}^j(k-d)&=&\Lambda_{l_{k-d},l_k}[B_{l_k}'(P_{l_k}(k)-P_{l_k}^0(k))F_{l_k,l_{k-d+j+1}}B_{l_{k-d+j}}]\nonumber\\
&&-\sum_{s=1}^{d-j}\{\Lambda_{l_{k-d},l_{k-s}}[(T_{l_{k-s}}^s(k-s))'W_{l_{k-s}}(k-s)^{-1}T_{l_{k-s}}^0(k-s)\nonumber\\
&&\times F_{l_{k-s},l_{k-d+j+1}}B_{l_{k-d+j}}]\}\nonumber\\
&&-\sum_{s=d-j+1}^{d-1}\{\Lambda_{l_{k-d},l_{k-s}}[(T_{l_{k-s}}^s(k-s))'W_{l_{k-s}}(k-s)^{-1} T_{l_{k-s}}^{s-(d-j)}(k-s)]\}, \nonumber\\
&&j=2,3,\cdots,d-1 \label{a4}
\end{eqnarray}
for $k=N, N-1,\cdots,0, l_{k-d}\in \Theta$ with terminal values
\begin{eqnarray}
&&T_{l_{N-i}}^j(N-i)=0,\nonumber\\
&& j=0,1,\cdots,d-1, i=0,1,\cdots,d-1, l_{N-i}\in \Theta,\nonumber
\end{eqnarray}
where $F_{l_k,l_i}=A_{l_k}\cdots A_{l_i} (i=0,\cdots,k), F_{l_k,l_{k+1}}=I$, and $P_{l_{k-1}}(k-1)$ and $P_{l_{k-1}}^0(k-1)$ satisfy the following backward recursions
\begin{eqnarray}
P_{l_{k-1}}(k-1)&=&\Lambda_{l_{k-1},l_k}[Q+A_{l_k}'(P_{l_k}(k)-P_{l_k}^0(k))A_{l_k}],\label{a5}\\
P_{l_{k-1}}^0(k-1)&=&(T_{l_{k-1}}^0(k-1))'W_{l_{k-1}}(k-1)^{-1}T_{l_{k-1}}^0(k-1)\label{a6}
\end{eqnarray}
for $k=N, N-1,\cdots,0, l_{k-1}\in \Theta$ with terminal values
\begin{eqnarray}
P_{l_N}(N)&=&\Lambda_{l_{N},l_{N+1}}P_{l_{N+1}}(N+1),\nonumber\\
P_{l_{N-i}}^0(N-i)&=&0, i=0,1,\cdots,d-1, l_{N-i}\in \Theta.\nonumber
\end{eqnarray}
(\ref{a1})-(\ref{a6}) is termed as the backward coupled difference Riccati equation(CDRE). In addition, we introduce the following notations
\begin{eqnarray}
(\alpha_{l_{k-1},l_{k-1}}^{d-1}(k-1,k-1))'&=&(\delta_{l_{k-1}}^{d-1}(k-1))',\label{a9}\\
(\alpha_{l_{k-1},l_{k-j}}^{d-j}(k-1,k-j))'&=&(\delta_{l_{k-1}}^{d-j}(k-1))'-\sum_{s=1}^{j-1}(\alpha_{l_{k-1},l_{k-s}}^{d-s}(k-1,k-s))'\nonumber\\
&&\times W_{l_{k-s-1}}(k-s-1)^{-1}T_{l_{k-s-1}}^{d-j+s}(k-s-1),\nonumber\\
&&j=2,3,\cdots,d-1\label{a10}\\
(\delta_{l_{k-1}}^1(k-1))'&=&\Lambda_{l_{k-1},l_k}[A_{l_k}'(P_{l_k}(k)-P_{l_k}^0(k))B_{l_k}]-(T_{l_{k-1}}^0(k-1))'\nonumber\\
&&\times W_{l_{k-1}}(k-1)^{-1}T_{l_{k-1}}^1(k-1),\label{a7}\\
(\delta_{l_{k-1}}^j(k-1))'&=&\Lambda_{l_{k-1},l_k}[A_{l_k}'(\delta_{l_k}^{j-1}(k))']-(T_{l_{k-1}}^0(k-1))'W_{l_{k-1}}(k-1)^{-1}\nonumber\\
&&\times T_{l_{k-1}}^j(k-1),\nonumber\\
&& j=2,3,\cdots,d-1, \label{a8}
\end{eqnarray}
for $k=N, N-1,\cdots,0, l_{k-1}, l_{k-j}\in \Theta$. The results for the finite horizon optimal control are stated as below.
\begin{lemma}
Consider $\alpha_{l_{k-1},l_{k-j}}^{d-j}(k-1,k-j)$ and $T_{l_{k-j}}^j(k-j)$ as in (\ref{a2})-(\ref{a4}), (\ref{a9}) and (\ref{a10}), the following expressions are satisfied
\begin{eqnarray}
\mbox{E}\{A_{l_k}'(\alpha_{l_k,l_k}^{d-1}(k,k))'|\mathcal{G}_{k-1}\}&=&(T_{l_{k-1}}^0(k-1))',\label{c1}\\
\mbox{E}\{A_{l_k}'(\alpha_{l_k,l_{k-j+1}}^{d-j}(k,k-j+1))'|\mathcal{G}_{k-1}\}&=&(\alpha_{l_{k-1},l_{k-j+1}}^{d-j+1}(k-1,k-j+1))', \nonumber\\
&&j=2,\cdots,d-1,\label{c2}\\
\mbox{E}\{B_{l_k}'(\alpha_{l_k,l_k}^{d-1}(k,k))'|\mathcal{G}_{k-1}\}&=&(T_{l_{k-1}}^1(k-1))',\label{c3}\\
\mbox{E}\{B_{l_k}'(\alpha_{l_k,l_{k-j+1}}^{d-j}(k,k-j+1))'|\mathcal{G}_{k-j}\}&=&(T_{l_{k-j}}^j(k-j))', j=2,\cdots,d-1.\label{c4}
\end{eqnarray}
\end{lemma}
\emph{Proof}. The detailed proof can be found in \cite{HLZ2018}.

\begin{theorem}
Under Assumption 1, the finite-horizon optimal controller is given by
\begin{eqnarray}
u(k-d)&=&-W_{l_{k-d}}(k-d)^{-1}T_{l_{k-d}}^0(k-d)x(k-d+1)\nonumber\\
&&-\sum_{j=1}^{d-1}W_{l_{k-d}}(k-d)^{-1}T_{l_{k-d}}^j(k-d) u(k-2d+j)\label{a13}\\
&&k=d,d+1,\cdots,N, l_{k-d}\in \Theta\nonumber
\end{eqnarray}
and the controller is unique if and only if
\begin{eqnarray}
W_{l_{k-d}}(k-d)>0, \ k=N,N-1,\cdots,d, l_{k-d}\in \Theta.\nonumber
\end{eqnarray}

The optimal costate is
\begin{eqnarray}
\lambda_{k-1}&=&(P_{l_{k-1}}(k-1)-P_{l_{k-1}}^0(k-1))x(k)-\sum_{s=1}^{d-1}(\alpha_{l_{k-1},l_{k-s}}^{d-s}(k-1,k-s))'\nonumber\\
&&\times W_{l_{k-s-1}}(k-s-1)^{-1}\mbox{E}\{\alpha_{l_{k-1},l_{k-s}}^{d-s}(k-1,k-s)x(k)|\mathcal{G}_{k-s-1}\}\label{a14}
\end{eqnarray}
and the optimal cost is
\begin{eqnarray}
J_N^{*}&=&\mbox{E}\{\sum_{k=0}^{d-1}x(k)'Qx(k)+x(d)'(P_{l_{d-1}}(d-1)-P_{l_{d-1}}^0(d-1))x(d)\nonumber\\
&&-x(d)'\sum_{s=1}^{d-1}(\alpha_{l_{d-1},l_{d-s}}^{d-s}(d-1,d-s))'W_{l_{d-s-1}}(d-s-1)^{-1}\mbox{E}[\alpha_{l_{d-1},l_{d-s}}^{d-s}(d-1,d-s)\nonumber\\
&&\times x(d)|\mathcal{G}_{d-s-1}]\}.\label{ff24}
\end{eqnarray}

\end{theorem}
Proof. The detailed proof can be found in \cite{HLZ2018}.

\begin{remark}
 The finite horizon optimal control result is obtained by two basic formulas: one is an improved delayed forward and backward jumping parameter equation (D-FBJPE) which is used to deal with the input delay, and the other is a d-step backward formula which is used to overcome the correlation of the jumping parameters. The detailed description can be seen from \cite{HLZ2018}.
\end{remark}

\begin{remark}
In order to make the time horizon $N$ explicit in the finite-time optimal control problem, we rewrite $W_{l_{k-d}}(k-d),T_{l_{k-d}}^i(k-d)(i=0,1,\cdots,d-1), P_{l_{k-1}}(k-1), P_{l_{k-1}}^0(k-1),\delta_{l_{k-1}}^j(k-1)(j=1,\cdots,d-1), \alpha_{l_{k-1},l_{k-j}}^{d-j}(k-1,k-j)(j=1,\cdots,d-1)$ as $W_{l_{k-d}}(k-d,N),T_{l_{k-d}}^i(k-d,N)(i=0,1,\cdots,d-1), P_{l_{k-1}}(k-1,N), P_{l_{k-1}}^0(k-1,N),\delta_{l_{k-1}}^j(k-1,N)(j=1,\cdots,d-1), \alpha_{l_{k-1},l_{k-j}}^{d-j}(k-1,k-j,N)(j=1,\cdots,d-1)$ respectively, and set $P_{l_{N+1}}(N+1)=0 (l_{N+1}=1,\cdots,N)$.
\end{remark}

\section{Main Results}

In this section, the main results of this paper will be presented, the necessary and sufficient stabilization conditions for optimal control systems will be established.

Before proposing the solution to Problem 1, the following lemmas will be given first.

\begin{lemma}
If $R>0$, the optimal control for the finite horizon case $(N\geq d)$ has a unique solution.
\end{lemma}
\emph{Proof}. See Appendix A.

\begin{remark}
If $R>0$, we will obtain from Lemma 2 that
\begin{eqnarray}
W_{l_{k-d}}(k-d,N+1)>0, k=N+1,\cdots,d.\nonumber
\end{eqnarray}
It follows from (\ref{a1})-(\ref{a4}) that $W_{l_{k-d}}(k-d,N+1)$ can be computed for $k=d-1,d-2,\cdots,0$ as well. Moreover, one yields from $P_{N+1}=0$ that
\begin{eqnarray}
W_{l_{k-1-d}}(k-1-d,N)=W_{l_{k-d}}(k-d,N+1)>0.\nonumber
\end{eqnarray}
So we have $W_{l_{k-d}}(k-d,N)>0$ for any $k=d-1,\cdots,0$ and $N$.
\end{remark}
\begin{lemma}
Take $N\geq d$. If $R>0$, we have
\begin{eqnarray}
&&P_{l_{k-1}}(k-1,N)\geq 0,k=N,\cdots,1, l_{k-1}\in \Theta,\label{s9}\\
&&P_{l_{k-1}}^0(k-1,N)\geq 0,k=N,\cdots,1, l_{k-1}\in \Theta,\label{s10}\\
&&(P_{l_{k-1}}(k-1,N)-P_{l_{k-1}}^0(k-1,N))-\sum_{s=1}^{d-1}(\alpha_{l_{k-1},l_{k-s}}^{d-s}(k-1,k-s,N))'\nonumber\\
&&\times W_{l_{k-1-s}}(k-1-s,N)^{-1}\alpha_{l_{k-1},l_{k-s}}^{d-s}(k-1,k-s,N)\geq 0,\label{s11}\\
&&k=N,\cdots,d, l_{k-1},l_{k-d}\in \Theta.\nonumber
\end{eqnarray}
\end{lemma}
\emph{Proof}. See Appendix B.

\begin{lemma}
If Assumptions 1 and 2 are satisfied, then there exists an integer $N_0\geq d$, such that
\begin{eqnarray}
&&(P_{l_{d-1}}(d-1,N_0)-P_{l_{d-1}}^0(d-1,N_0))-\sum_{s=1}^{d-1}(\alpha_{l_{d-1},l_{d-s}}^{d-s}(d-1,d-s,N_0))'\nonumber\\
&&\times W_{l_{d-1-s}}(d-1-s,N_0)^{-1}\alpha_{l_{d-1},l_{d-s}}^{d-s}(d-1,d-s,N_0)>0,\label{s15}\\
&&l_{d-s}\in \Theta, s=1,\cdots,d-1.\nonumber
\end{eqnarray}
\end{lemma}
\emph{Proof}. See Appendix C.

In what follows, we introduce the generalized CAREs as
\begin{eqnarray}
W_{l_d}&=&\Lambda_{l_d,l_0}[B_{l_0}'(P_{l_0}-P_{l_0}^0)B_{l_0}+R_{l_0}] -\sum_{s=1}^{d-1}\{\Lambda_{l_d,l_s}[(T_{l_{s}}^s)'W_{l_{s}}^{-1}
(T_{l_{s}}^s)]\},\label{s19}\\
T_{l_{d}}^0&=&\Lambda_{l_d,l_0}[B_{l_0}'(P_{l_0}-P_{l_0}^0)A_{l_0}A_{l_{1}}\cdots A_{l_{d-1}}]\nonumber\\
&&-\sum_{s=1}^{d-1}\{\Lambda_{l_d,l_s}[(T_{l_{s}}^s)'(W_{l_{s}})^{-1}T_{l_{s}}^0 A_{l_{s}}\cdots A_{l_{d-1}}]\},\label{s20}\\
T_{l_{d}}^1&=&\Lambda_{l_d,l_0}[B_{l_0}'(P_{l_0}-P_{l_0}^0)A_{l_0}\cdots A_{l_{d-2}}B_{l_{d-1}}]\nonumber\\
&&-\sum_{s=1}^{d-2}\{\Lambda_{l_d,l_s}[(T_{l_{s}}^s)'(W_{l_{s}})^{-1}T_{l_{s}}^0A_{l_{s}}\cdots A_{l_{d-2}}B_{l_{d-1}}]\}\nonumber\\
&&-\sum_{l_{d-1}^{d-1}}\{\Lambda_{l_d,l_{d-1}}[(T_{l_{d-1}}^{d-1})'(W_{l_{d-1}})^{-1}T_{l_{d-1}}^0B_{l_{d-1}}]\},\label{s21}\\
T_{l_{d}}^j&=&\Lambda_{l_d,l_0}[B_{l_0}'(P_{l_0}-P_{l_0}^0)A_{l_0}\cdots A_{l_{d-j-1}}B_{l_{d-j}}]\nonumber\\
&&-\sum_{s=1}^{d-j-1}\{\Lambda_{l_d,l_s}[(T_{l_{s}}^s)'(W_{l_{s}})^{-1}T_{l_{s}}^0A_{l_{s}}\cdots A_{l_{d-j-1}}B_{l_{d-j}}]\}\nonumber\\
&&-\Lambda_{l_d,l_{d-j}}[(T_{l_{d-j}}^{d-j})'(W_{l_{d-j}})^{-1}T_{l_{d-j}}^0B_{l_{d-j}}]\nonumber\\
&&-\sum_{s=d-j+1}^{d-1}\{\Lambda_{l_d,l_s}[(T_{l_{s}}^s)'(W_{l_{s}})^{-1}T_{l_{s}}^{s-(d-j)}]\}, j=2,3,\cdots,d-1, \label{s22}
\end{eqnarray}
\begin{eqnarray}
P_{l_{1}}&=&\Lambda_{l_1,l_0}[Q_{l_0}+A_{l_0}'(P_{l_0}-P_{l_0}^0)A_{l_0}],\label{s23}\\
P_{l_{1}}^0&=&(T_{l_{1}}^0)'(W_{l_{1}})^{-1}T_{l_{1}}^0,\label{s24}\\
(\delta_{l_{1}}^1)'&=&\Lambda_{l_1,l_0}[A_{l_0}'(P_{l_0}-P_{l_0}^0)B_{l_0}]-(T_{l_{1}}^0)'(W_{l_{1}})^{-1}T_{l_{1}}^1,\label{s25}\\
(\delta_{l_{1}}^j)'&=&\Lambda_{l_1,l_0}[A_{l_0}'(\delta_{l_0}^{j-1})']-(T_{l_{1}}^0)'(W_{l_{1}})^{-1}T_{l_{1}}^j,j=2,3,\cdots,d-1, \label{s26}\\
(\alpha_{l_{1},l_{1}}^{d-1})'&=&(\delta_{l_{1}}^{d-1})',\label{s27}\\
(\alpha_{l_{1},l_{j}}^{d-j})'&=&(\delta_{l_{1}}^{d-j})'-\sum_{s=1}^{j-1}(\alpha_{l_{1},l_{s}}^{d-s})'(W_{l_{s+1}})^{-1}T_{l_{s+1}}^{d-j+s},
j=2,3,\cdots,d-1.\label{s28}
\end{eqnarray}

Then the main results of this paper can be stated as below.
\begin{theorem}
If Assumptions 1 and 2 are satisfied and system (\ref{fff1}) is mean square stabilizable, we have the following properties:

1) When $N\rightarrow \infty$, $P_{l_{k-1}}(k-1,N), P_{l_{k-1}}^0(k-1,N), W_{l_{k-d}}(k-d,N), T_{l_{k-d}}^s(k-d,N) (s=0,1,\cdots,d-1)$ converge to $P_{l_1}, P_{l_1}^0, W_{l_d}, T_{l_d}^s(s=0,1,\cdots,d-1)$ respectively for any $k\geq 0$ and $l_{k-1}\in \Theta$. Furthermore, $P_{l_1}, P_{l_1}^0, W_{l_d}, T_{l_d}^s(s=0,1,\cdots,d-1)$ obey the coupled algebraic Riccati equations (\ref{s19})-(\ref{s28}).

2)
\begin{eqnarray}
(P_{l_1}-P_{l_1}^0)-\sum_{s=1}^{d-1}(\alpha_{l_1,l_s}^{d-s})'W_{l_{s+1}}^{-1}\alpha_{l_1,l_s}^{d-s}>0.\label{ss29}
\end{eqnarray}
\end{theorem}
\emph{Proof}. See Appendix D.

\begin{theorem}
If Assumptions 1 and 2 are satisfied, then system (\ref{fff1}) is mean-square stabilizable if and only if there exists a unique solution to (\ref{s19})-(\ref{s28}) such that
\begin{eqnarray}
(P_{l_1}-P_{l_1}^0)-\sum_{s=1}^{d-1}(\alpha_{l_1,l_s}^{d-s})'W_{l_{s+1}}^{-1}\alpha_{l_1,l_s}^{d-s}>0.\nonumber
\end{eqnarray}
In this case, the optimal controller given by
\begin{eqnarray}
u(k)=-W_{l_d}^{-1}T_{l_d}^0A_{l_d}x(k)-W_{l_d}^{-1}T_{l_d}^0B_{l_d}u(k-d)-\sum_{j=1}^{d-1}W_{l_d}^{-1}T_{l_d}^ju(k-d+j),\label{s57}
\end{eqnarray}
stabilizes (\ref{fff1}) and minimizes the performance index (\ref{ss1}). The optimal value of (\ref{ss1}) is given by
\begin{eqnarray}
J^*&=&\mbox{E}\{x_0'(P_{l_1}-P_{l_1}^0)x_0-x_0'\sum_{s=1}^{d-1}(\alpha_{l_1,l_s}^{d-s})'W_{l_{s+1}}^{-1}\mbox{E}[\alpha_{l_1,l_s}^{d-s}x_0|{\cal{G}}_{-s-1}]\nonumber\\
&&-\sum_{k=0}^{d-1}u(k-d)'Ru(k-d)+\sum_{k=0}^{d-1}[u(k-d)+W_{l_d}^{-1}T_{l_d}^0x(k-d+1)\nonumber\\
&&+\sum_{j=1}^{d-1}W_{l_d}^{-1}T_{l_d}^ju(k-2d+j)]'W_{l_d}^{-1}[u(k-d)+W_{l_d}^{-1}T_{l_d}^0x(k-d+1)\nonumber\\
&&+\sum_{j=1}^{d-1}W_{l_d}^{-1}T_{l_d}^ju(k-2d+j)]\}\label{s58}
\end{eqnarray}
\end{theorem}
\emph{Proof}. See Appendix E.

\section{Numerical Examples}
In this section, we show the efficiency of the stabilization result for the infinite-horizon case. The specifications of the system (\ref{fff1}) and the weighting matrices in (\ref{ss1}) are as follows
\begin{eqnarray}
&&A_1=\left[
      \begin{array}{cc}
        2 & 1.1 \\
        -1.7 & -0.8 \\
      \end{array}
    \right], A_2=\left[
                    \begin{array}{cc}
                      0.8 & 0 \\
                      0 & 0.6 \\
                    \end{array}
                  \right], B_1=\left[
                                  \begin{array}{c}
                                    1 \\
                                    1 \\
                                  \end{array}
                                \right],B_2=\left[
                                               \begin{array}{c}
                                                 2 \\
                                                 1 \\
                                               \end{array}
                                             \right],\nonumber\\
&&Q=\left[
          \begin{array}{cc}
            1& 0 \\
            0 & 1 \\
          \end{array}
        \right],R=1,P_{N+1}=\left[
          \begin{array}{cc}
            1& 0 \\
            0 & 1 \\
          \end{array}
        \right].\nonumber
\end{eqnarray}

The initial distribution of $\theta(k)$ is $(0.5,~0.5)$s and probability transition matrix of $\theta(k)$ is $\left[
                                      \begin{array}{cc}
                                        0.9 & 0.1 \\
                                        0.3 & 0.7 \\
                                      \end{array}
                                    \right]$. Meanwhile, the constant delay $d=2$, and the initial values of $x_0,u(-1),u(-2)$ remain unchanged. we run $50$
Monte Carlo simulations, and select the first trajectory to show the efficiency of the proposed algorithm. By applying Theorem 2, the calculation results are listed as follow.
\begin{eqnarray}
P_1&=&\left[
        \begin{array}{cc}
          23.2324 & 13.0039 \\
          13.0039 & 8.7742 \\
        \end{array}
      \right],P_2=\left[
        \begin{array}{cc}
          12.0410 & 5.5718 \\
          5.5718 & 4.4655 \\
        \end{array}
      \right],\nonumber\\
 P_1^0&=&\left[
        \begin{array}{cc}
          6.4216 & 3.9763 \\
          3.9763 & 2.4622 \\
        \end{array}
      \right],P_2^0=\left[
        \begin{array}{cc}
          3.5328 & 1.7056 \\
          1.7056 & 0.8235 \\
        \end{array}
      \right],\nonumber \\
      W_1&=&23.4331, W_2=26.0873\nonumber\\
      T_1^0&=&\left[
                \begin{array}{cc}
                  12.2669 & 7.5958 \\
                \end{array}
              \right], T_2^0=\left[
                \begin{array}{cc}
                  9.6 & 4.6349 \\
                \end{array}
              \right],\nonumber\\
      T_1^1&=&21.8551, T_2^1=24.5677.\nonumber
\end{eqnarray}
And the stabilization condition
\begin{eqnarray}
P_1-P_1^0-(\alpha_{1,1}^1)'W_1^{-1}\alpha_{1,1}^1&=&\left[
                                                      \begin{array}{cc}
                                                        9.1523 & 4.3331 \\
                                                        4.3331 & 3.4345 \\
                                                      \end{array}
                                                    \right]>0,\nonumber\\
P_1-P_1^0-(\alpha_{1,1}^1)'W_2^{-1}\alpha_{1,1}^1&=&\left[
                                                      \begin{array}{cc}
                                                        10.3881 & 5.0906 \\
                                                        5.0906 & 3.8988 \\
                                                      \end{array}
                                                    \right]>0,\nonumber\\
P_2-P_2^0-(\alpha_{2,2}^1)'W_1^{-1}\alpha_{2,2}^1&=&\left[
                                                      \begin{array}{cc}
                                                        3.8639 & 1.5032 \\
                                                        1.5032 & 2.4397 \\
                                                      \end{array}
                                                    \right]>0,\nonumber\\
P_2-P_2^0-(\alpha_{2,2}^1)'W_2^{-1}\alpha_{2,2}^1&=&\left[
                                                      \begin{array}{cc}
                                                        4.6133 & 1.8845 \\
                                                        1.8845 & 2.6337 \\
                                                      \end{array}
                                                    \right]>0\nonumber
\end{eqnarray}
is satisfied. If $\theta(k)=1$, the optimal stabilization controller is
\begin{eqnarray}
u(k)
&=&-W_1^{-1}T_1^0x(k+1)-W_1^{-1}T_1^1u(k-1)\nonumber\\
&=&-\left[
                                             \begin{array}{cc}
                                               0.5235 & 0.3242 \\
                                             \end{array}
                                           \right]x(k+1)-0.9327u(k-1).\nonumber
\end{eqnarray}
If $\theta(k)=2$, the optimal stabilization controller is
\begin{eqnarray}
u(k)&=&-W_2^{-1}T_2^0x(k+1)-W_2^{-1}T_2^1u(k-1)\nonumber\\
&=&-\left[
                                             \begin{array}{cc}
                                               0.3680 & 0.1777 \\
                                             \end{array}
                                           \right]x(k+1)-0.9417u(k-1).\nonumber
\end{eqnarray}
From (\ref{s58}), we can derive the optimal cost $J_0=287.4952$.
%

%
\section{Conclusion}
In this paper, necessary and sufficient stabilization conditions
for MJLS with input delay have been investigated. It is shown that,  under the exact observability assumption, we
show that the closed MJLS is mean square stabilizable if
and only if the CARE admits a unique positive definite
solution.


\appendix

\section{Proof of Lemma 2}
\emph{Proof}. By induction method, we will show that the optimal control associated with performance (\ref{ss1}) and system (\ref{fff1}) has a unique solution for $N\geq d$. For $N=d$,
\begin{eqnarray}
W_{l_{0}}(0,N)&=&\Lambda_{l_0,l_N}[B_{l_N}'(P_{l_N}(N,N)-P_{l_N}^0(N,N))B_{l_N}+R]\nonumber\\
&&-\sum_{s=1}^{N-1}\{\Lambda_{l_0,l_{N-s}}[(T_{l_{N-s}}^s(N-s,N))'W_{l_{N-s}}(N-s,N)^{-1}\nonumber\\
&&\times T_{l_{N-s}}^s(N-s,N)]\}.\nonumber
\end{eqnarray}
It follows from
\begin{eqnarray}
P_{l_N}(N,N)&=&\Lambda_{l_N,l_{N+1}}P_{l_{N+1}}(N+1)=0,\nonumber\\
P_{l_N}^0(N,N)&=&0,\nonumber\\
T_{l_{N-i}}^j(N-i,N)&=&0,i,j=0,1,\cdots,N-1\nonumber
\end{eqnarray}
that
\begin{eqnarray}
W_{l_0}(0,N)=R>0.\nonumber
\end{eqnarray}
In light of Theorem 1, we know that the finite horizon optimal control admits a unique solution for $N=d$.

Next, suppose the solution to finite horizon control with terminal time $M=m$ is unique for some $m\geq d$, i.e.,
\begin{eqnarray}
W_{l_{k-d}}(k-d,m)>0,k=m,\cdots,d.\label{s3}
\end{eqnarray}
Without loss of generality, we take $x(0)=0,u(-2)=u(-3)=\cdots=u(-d)=0$ and $u(-1)$ is an arbitrary value. Then we get that $x(1)=\cdots=x(d-1)=0$ and $x(d)=B_{\theta(d-1)}u(-1)$ is arbitrary. In view of (\ref{ff24}), we can obtain the optimal cost with terminal time $m$
\begin{eqnarray}
J_m^*
&=&u(-1)'\{\Lambda_{l_{-1},l_{d-1}}[B_{l_{d-1}}'(P_{l_{d-1}}(d-1,m)-P_{l_{d-1}}^0(d-1,m))B_{l_{d-1}}]\nonumber\\
&&-\sum_{s=1}^{d-1}\{\Lambda_{l_{-1},l_{d-1-s}}[T_{l_{d-1-s}}^s(d-1-s,m)'W_{l_{d-1-s}}(d-1-s,m)^{-1}\nonumber\\
&&\times T_{l_{d-1-s}}^s(d-1-s,m)]\}\}u(-1)\geq 0,\label{s5}
\end{eqnarray}
Since $u(-1)$ is arbitrary, we have
\begin{eqnarray}
&&\Lambda_{l_{-1},l_{d-1}}[B_{l_{d-1}}'(P_{l_{d-1}}(d-1,m)-P_{l_{d-1}}^0(d-1,m))B_{l_{d-1}}]\nonumber\\
&&\geq\sum_{s=1}^{d-1}\{\Lambda_{l_{-1},l_{d-1-s}}[T_{l_{d-1-s}}^s(d-1-s,m)'W_{l_{d-1-s}}(d-1-s,m)^{-1}\nonumber\\
&&\times T_{l_{d-1-s}}^s(d-1-s,m)]\}\geq 0.\label{s7}
\end{eqnarray}
Recalling that the variables defined in (\ref{a1})-(\ref{a4}) are time-invariant for $N$ owning to the selection that $P_{N+1}=0$, i.e.,
\begin{eqnarray}
W_{l_{k-d}}(k-d,m)&=&W_{l_{k-d-s}}(k-d-s,m-s), \nonumber\\
T_{l_{k-d}}^0(k-d,m)&=&T_{l_{k-d-s}}^0(k-d-s,m-s),\nonumber\\
&\cdots&\nonumber\\
T_{l_{k-d}}^{d-1}(k-d,m)&=&T_{l_{k-d-s}}^{d-1}(k-d-s,m-s),\nonumber\\
P_{l_{k-1}}(k-1,m)&=&P_{l_{k-1-s}}(k-1-s,m-s),\nonumber\\
P_{l_{k-1}}^0(k-1,m)&=&P_{l_{k-1-s}}^0(k-1-s,m-s).\nonumber
\end{eqnarray}
For $N=m$, it yields from (\ref{s3}) that
\begin{eqnarray}
W_{l_{k-d}}(k-d,m+1)=W_{l_{k-d-1}}(k-d-1,m)>0,k=m+1,\cdots,d+1.\label{s8}
\end{eqnarray}
For $k=d$, we have
\begin{eqnarray}
W_{l_0}(0,m+1)&=&\Lambda_{l_0,l_d}[B_{l_d}'(P_{l_d}(d,m+1)-P_{l_d}^0(d,m+1))B_{l_d}+R]\nonumber\\
&&-\sum_{s=1}^{d-1}\{\Lambda_{l_0,l_{d-s}}(T_{l_{d-s}}^s(d-s,m+1))'W_{l_{d-s}}(d-s,m+1)^{-1}\nonumber\\
&&\times T_{l_{d-s}}^s(d-s,m+1)\}\nonumber\\
&=&\Lambda_{l_{-1},l_{d-1}}[B_{l_{d-1}}'(P_{l_{d-1}}(d-1,m)-P_{l_{d-1}}^0(d-1,m))B_{l_{d-1}}+R]\nonumber\\
&&-\sum_{s=1}^{d-1}\{\Lambda_{l_{-1},l_{d-1-s}}[(T_{l_{d-1-s}}^s(d-1-s,m))'\nonumber\\
&&\times W_{l_{d-1-s}}(d-1-s,m)^{-1}T_{l_{d-1-s}}^s(d-1-s,m)]\}\nonumber\\
&\geq&0.\nonumber
\end{eqnarray}

It concludes from Theorem 1 that there exists a unique solution to the finite horizon optimal control with $N=m+1$. Then the uniqueness and existence of the solution to the optimal control for the finite horizon case is shown for $N\geq d$. This completes the proof.

\section{Proof of Lemma 3}

\emph{Proof}. In view of Theorem 1, we have that $W_{l_{k-d}}(k-d)>0$ for any $N\geq d$ and $d\leq k\leq N$. Therefore, we have from (\ref{a6}) that $P_{l_{k-1}}^0(k-1,N)\geq 0$ for $k=N,\cdots,1, l_{k-1}\in \Theta$. Recalling from the invariance of $P_{l_{k-1}}(k-1,N)=P_{l_{d-1}}(d-1,N-k+d)$
and (\ref{s7}), one yields
\begin{eqnarray}
\Lambda_{l_{k-d-1,l_{k-1}}}[B_{l_{k-1}}'(P_{l_{k-1}}(k-1,N)-P_{l_{k-1}}^0(k-1,N))B_{l_{k-1}}]\geq 0.\nonumber
\end{eqnarray}
It implies that $P_{l_{k-1}}(k-1,N)-P_{l_{k-1}}^0(k-1,N)\geq 0$. Then $P_{l_{k-1}}(k-1,N)\geq P_{l_{k-1}}^0(k-1,N)\geq 0$ follows immediately.

Next, we will show that (\ref{s11}) is satisfied. Let the system (\ref{fff1}) start at $d$ with any initial value $x_d$ and denote it as
\begin{eqnarray}
S_d=\sum_{k=d}^N\mbox{E}\{x(k)'Qx(k)+u(k-d)'Ru(k-d)\}.\label{s12}
\end{eqnarray}
In view of (\ref{ff24}), the optimal value of (\ref{s12}) can be written as
\begin{eqnarray}
S_d^*
&=&x(d)'\{(P_{l_{d-1}}(d-1,N)-P_{l_{d-1}}^0(d-1,N))\nonumber\\
&&-\sum_{s=1}^{d-1}(\alpha_{l_{d-1},l_{d-s}}^{d-s}(d-1,d-s,N))'W_{l_{d-s-1}}(d-s-1,N)^{-1}\nonumber\\
&&\times \alpha_{l_{d-1},l_{d-s}}^{d-s}(d-1,d-s,N)\}x(d)\geq 0.\label{s13}
\end{eqnarray}
It follows from the arbitrary of $x(d)$ that
\begin{eqnarray}
&&(P_{l_{d-1}}(d-1,N)-P_{l_{d-1}}^0(d-1,N))\nonumber\\
&&-\sum_{s=1}^{d-1}(\alpha_{l_{d-1},l_{d-s}}^{d-s}(d-1,d-s,N))'W_{l_{d-s-1}}(d-s-1,N)^{-1}\nonumber\\
&&\times \alpha_{l_{d-1},l_{d-s}}^{d-s}(d-1,d-s,N)\geq 0.\label{s14}
\end{eqnarray}
Based on (\ref{s14}) and the time invariance of $P_{l_{k-1}}(k-1,N), P_{l_{k-1}}^0(k-1,N), \alpha_{l_{k-1},l_{k-s}}^{d-s}(k-1,k-s,N), W_{l_{k-s-1}}(k-s-1,N)$, one yields (\ref{s11}). This completes the proof of Lemma 3.

\section{Proof of Lemma 4}

\emph{Proof}. In order to facilitate the description, we denote
\begin{eqnarray}
&&\Omega_{l_{d-1},l_{1}}(d-1,1,N)\nonumber\\
&=&(P_{l_{d-1}}(d-1,N)-P_{l_{d-1}}^0(d-1,N))-\sum_{s=1}^{d-1}(\alpha_{l_{d-1},l_{d-s}}^{d-s}(d-1,d-s,N))'\nonumber\\
&&\times W_{l_{d-1-s}}(d-1-s,N)^{-1}\alpha_{l_{d-1},l_{d-s}}^{d-s}(d-1,d-s,N).\nonumber
\end{eqnarray}
If Assumption 1 is satisfied, one yields from Lemma 3 that $\Omega_{l_{d-1},l_{1}}(d-1,1,N)\geq0$ for all $N\geq d$. In the next, we just need to show that there exists $N_0\geq d$ such that $\Omega_{l_{d-1},l_{1}}(d-1,1,N_0)>0$. Assume this is not valid. Then we get an non-empty set
\begin{eqnarray}
X_N\triangleq \{x\in R^n: x\neq 0,x'\Omega_{l_{d-1},l_{1}}(d-1,1,N)x=0\}.\nonumber
\end{eqnarray}
In light of (\ref{s12}) and (\ref{s13}), we can deduce that $x'\Omega_{l_{d-1},l_{1}}(d-1,1,N)x\leq x'\Omega_{l_{d-1},l_{1}}(d-1,1,N+1)x$. Since $x$ is arbitrary, we get that $\Omega_{l_{d-1},l_{1}}(d-1,1,N)\leq \Omega_{l_{d-1},l_{1}}(d-1,1,N+1)$. Then if $x'\Omega_{l_{d-1},l_{1}}(d-1,1,N+1)x=0$, we can deduce that $x'\Omega_{l_{d-1},l_{1}}(d-1,1,N)x=0$. It implies that $X_{N+1}\subset X_N$.
Note that each $X_N$ is non-empty and with finite-dimension, we can obtain that
\begin{eqnarray}
1\leq \cdots\leq \dim(X_{d+2})\leq \dim(X_{d+1})\leq \dim(X_d)\leq n.\label{s16}
\end{eqnarray}
It follows from (\ref{s16}) that there must exist an integer $N_1$, such that for $N\geq N_1$, $\dim(X_N)=\dim(X_{N_1})$ and thus $X_N=X_{N_1}$. It means that $\bigcap_{N\geq d}X_N=X_{N_1}\neq \emptyset$. Therefore, there must exist a nonzero vector $x\in X_{N_1}$, such that
$x'\Omega_{l_{d-1},l_{1}}(d-1,1,N+1)x=0$ for any $N\geq d$.

Set $x(d)=x$ in (\ref{s13}), we get
\begin{eqnarray}
S_d^*&=&\min\sum_{k=d}^N\mbox{E}\{x(k)'Qx(k)+u(k-d)'Ru(k-d)\}\nonumber\\
&=&0.\label{s17}
\end{eqnarray}
It follows from the hypothesis $R>0$ and $Q=C'C\geq 0$ that
\begin{eqnarray}
u^*(k-d),Cx^*(k)=0,d\leq k\leq N,N\geq d.\nonumber
\end{eqnarray}
Then system (\ref{fff1}) becomes as
\begin{eqnarray}
x^*(k+1)&=&A_{\theta(k)}x^*(k),Cx^*(k)=0,\forall k\geq d.\label{s18}
\end{eqnarray}
From the observability of (\ref{s18}), we get that $x(d)=0$. This contradicts the fact $x\neq 0$. So there exists some $N_0\geq d$ such that $\Omega_{l_{d-1},l_{d-s}}(d-1,d-s,N_0)>0$. This completes the proof of Lemma 4.

\section{Proof of Theorem 2}

\emph{Proof}. (1) In the first part of the proof, we will show the convergence of the difference Riccati equations (\ref{a1})-(\ref{a4}).

Now, we start to prove that $W_{l_{k-d}}(k-d,N), T_{l_{k-d}}^s(k-d,N) (s=0,1,\cdots,d-1)$ are convergent. Define
\begin{eqnarray}
\bar{x}(k)=\mbox{col}\{x(k),u(k-1),\cdots,u(k-d)\},\nonumber
\end{eqnarray}
then the system (\ref{fff1}) can be rewritten as
\begin{eqnarray}
\bar{x}(k+1)=\bar{A}_{\theta(k)}\bar{x}(k)+\bar{B}u(k),\label{s29}
\end{eqnarray}
where
\begin{eqnarray}
\bar{A}_{\theta(k)}=\left[
                      \begin{array}{ccccc}
                        A_{\theta(k)} & 0 & \cdots & 0 & B_{\theta(k)} \\
                        0 & 0 & \cdots & 0 & 0 \\
                        0 & I & \cdots & 0 & 0 \\
                        \vdots & \vdots & \ddots & \vdots & \vdots \\
                        0 & 0 & \cdots & I & 0 \\
                      \end{array}
                    \right],\bar{B}=\left[
                                      \begin{array}{c}
                                        0 \\
                                        I \\
                                        0 \\
                                        \vdots \\
                                        0 \\
                                      \end{array}
                                    \right],\nonumber
\end{eqnarray}
and the performance index (\ref{ss1}) becomes as
\begin{eqnarray}
\bar{J}=\mbox{E}\{\sum_{k=0}^\infty\bar{x}(k)'\bar{Q}\bar{x}(k)+u(k)'Ru(k)\},\label{s30}
\end{eqnarray}
where
\begin{eqnarray}
\bar{Q}=\mbox{diag}\{Q,0,0,0,0\}.\nonumber
\end{eqnarray}
It can be seen that Problem 1 is equivalent to the minimization of (\ref{s30}) subject to (\ref{s29}). The necessary condition for minimizing the cost index of (\ref{s30}), i.e., the maximum principle, is given as
\begin{eqnarray}
0&=&\mbox{E}[\bar{B}\bar{\lambda}_k+Ru(k)|{\cal{G}}_k],\label{s31}\\
{\bar{\lambda}}_{k-1}&=&\mbox{E}[\bar{A}_{\theta(k)}'\bar{\lambda}_k+\bar{Q}_{\theta(k)}x(k)|{\cal{G}}_{k-1}],\label{s32}\\
\bar{\lambda}_N&=&\mbox{E}[\bar{P}_{\theta(N+1)}x(N+1)|{\cal{G}}_N].\label{s33}
\end{eqnarray}
Applying (\ref{s31})-(\ref{s33}) and following a similar derivation as that of Theorem 1 in \cite{ZLXF15}, we obtain the optimal controller
\begin{eqnarray}
u(k)=-\bar{\Upsilon}_{l_k}(k)^{-1}\bar{M}_{l_k}(k)\bar{x}(k),\label{s34}
\end{eqnarray}
and establish the relationship between $\bar{x}(k)$ and the costate $\bar{\lambda}_k$
\begin{eqnarray}
\bar{\lambda}_{k-1}=(\Lambda_{l_{k-1},l_k}\bar{P}_{l_k}(k))\bar{x}(k),\label{s35}
\end{eqnarray}
where $\bar{\Upsilon}_{l_k}(k), \bar{M}_{l_k}(k)$ and $\bar{P}_{l_k}(k)$ satisfy the following difference Riccati equations
\begin{eqnarray}
\bar{\Upsilon}_{l_{k}}(k)&=&\bar{B}'(\Lambda_{l_k,l_{k+1}}\bar{P}_{l_{k+1}}(k+1))\bar{B},\label{s36}\\
\bar{M}_{l_k}(k)&=&\bar{B}'(\Lambda_{l_k,l_{k+1}}\bar{P}_{l_{k+1}}(k+1))\bar{A}_{l_k},\label{s37}\\
\bar{P}_{l_k}(k)&=&\bar{A}_{l_k}'(\Lambda_{l_k,l_{k+1}}\bar{P}_{l_{k+1}}(k+1))\bar{A}_{l_k}+\bar{Q}-\bar{M}_{l_k}(k)'
\bar{\Upsilon}_{l_{k}}(k)\bar{M}_{l_k}(k).\label{s38}
\end{eqnarray}
According to the partitioned form of the augmented state $\bar{x}(k)$, the block forms of
$\bar{\lambda}_{k-1}$ and $\bar{P}_{l_k}(k)$ can be written as
\begin{eqnarray}
\bar{\lambda}_{k-1}&=&\mbox{col}\{\bar{\lambda}_{k-1}^{(0)},\bar{\lambda}_{k-1}^{(1)},\cdots,\bar{\lambda}_{k-1}^{(d)}\},\nonumber\\
\bar{P}_{l_k}(k)&=&\left[
                     \begin{array}{cccc}
                       \bar{P}_{l_k}^{(0,0)}(k) & \bar{P}_{l_k}^{(0,1)}(k) & \cdots & \bar{P}_{l_k}^{(0,d)}(k) \\
                       \bar{P}_{l_k}^{(1,0)}(k) & \bar{P}_{l_k}^{(1,1)}(k) & \cdots & \bar{P}_{l_k}^{(1,d)}(k) \\
                       \vdots & \vdots & \ddots & \vdots \\
                       \bar{P}_{l_k}^{(d,0)}(k) & \bar{P}_{l_k}^{(d,1)}(k) & \vdots & \bar{P}_{l_k}^{(d,d)}(k) \\
                     \end{array}
                   \right]
\end{eqnarray}
and (\ref{s35}) becomes as
\begin{eqnarray}
\bar{\lambda}_{k-1}^{(0)}&=&\Lambda_{l_{k-1},l_k}\{\bar{P}_{l_k}^{(0,0)}(k)x(k)+\bar{P}_{l_k}^{(0,1)}(k)u(k-1)+\cdots+\bar{P}_{l_k}^{(0,d)}(k)u(k-d)\},\label{s39}\\
\bar{\lambda}_{k-1}^{(1)}&=&\Lambda_{l_{k-1},l_k}\{\bar{P}_{l_k}^{(1,0)}(k)x(k)+\bar{P}_{l_k}^{(1,1)}(k)u(k-1)+\cdots+\bar{P}_{l_k}^{(1,d)}(k)u(k-d)\},\label{s40}\\
&\vdots&\nonumber\\
 \bar{\lambda}_{k-1}^{(d)}&=&\Lambda_{l_{k-1},l_k}\{\bar{P}_{l_k}^{(d,0)}(k)x(k)+\bar{P}_{l_k}^{(d,1)}(k)u(k-1)+\cdots+\bar{P}_{l_k}^{(d,d)}(k)u(k-d)\}.\label{s41}
\end{eqnarray}
From (\ref{s31}), one yields
\begin{eqnarray}
0=\mbox{E}\{\bar{\lambda}_{k}^{(1)}+Ru(k)|{\cal{G}}_{k}\}.\label{s42}
\end{eqnarray}
In addition, from the maximum principle developed in \cite{HLZ2018}, a necessary condition for minimizing (\ref{ss1}) for
system (\ref{fff1}) is as:
\begin{eqnarray}
0=\mbox{E}\{B_{\theta(k+d)}'\lambda_{k+d}+Ru(k)|{\cal{G}}_{k}\}.\label{s43}
\end{eqnarray}
Compared (\ref{s42}) with (\ref{s43}), we get
\begin{eqnarray}
\bar{\lambda}_k^{(1)}=\mbox{E}\{B_{\theta(k+d)}'\lambda_{k+d}|{\cal{G}}_{k}\},\nonumber
\end{eqnarray}
Based on the above expression, one yields
\begin{eqnarray}
T_{l_{k-1}}^0(k-1,N)&=&\Lambda_{l_{k-1},l_k}\bar{P}_{l_k}^{(1,0)}(k),\label{s49}\\
W_{l_{k-1}}(k-1,N)-R&=&\Lambda_{l_{k-1},l_k}\bar{P}_{l_k}^{(1,1)}(k),\label{s50}\\
T_{l_{k-1}}^{d-1}(k-1,N)&=&\Lambda_{l_{k-1},l_k}\bar{P}_{l_k}^{(1,2)}(k),\label{s51}\\
&&\vdots\nonumber\\
T_{l_{k-1}}^{1}(k-1,N)&=&\Lambda_{l_{k-1},l_k}\bar{P}_{l_k}^{(1,d)}(k).\label{s52}
\end{eqnarray}
Follow a similar discussion as (73)-(74) in \cite{ZLXF15}, we can show that if system (\ref{fff1}) is mean-square stabilizable, then the augmented system (\ref{s29}) is stabilizable in the mean square sense as well. Also, the observable keeps. So it follows from Theorem 2 in \cite{HLZ18} that $\bar{P}_{l_k}^{(1,0)}(k),$ $\bar{P}_{l_k}^{(1,1)}(k),$ $\cdots, \bar{P}_{l_k}^{(1,d)}(k)$ are convergent. It concludes from the ergodicity of $\theta(k)$ and (\ref{s49})-(\ref{s52}) that $T_{l_{k-1}}^0(k-1,N),\cdots,T_{l_{k-1}}^{d-1}(k-1,N), W_{l_{k-1}}(k-1,N)-R$ converge as well, i.e. $T_{l_{k-d}}^0(k-d,N-d+1),\cdots,T_{l_{k-d}}^{d-1}(k-d,N-d+1), W_{l_{k-d}}(k-d,N-d+1)-R$ are convergent. Denote $T_{l_d}^0=\lim_{N\rightarrow \infty}T_{l_{k-d}}^0(k-d,N-d+1),\cdots,T_{l_d}^{d-1}=\lim_{N\rightarrow \infty}T_{l_{k-d}}^{d-1}(k-d,N-d+1), W_{l_d}=\lim_{N\rightarrow \infty}W_{l_{k-d}}(k-d,N-d+1)$. Obviously, $T_{l_d}^0,\cdots,T_{l_d}^{d-1}, W_{l_d}$ satisfy (\ref{s19})-(\ref{s22}).

Next, we will show the convergence of $P_{l_{k-1}}(k-1,N)$ and $P_{l_{k-1}}^0(k-1,N)$. In view of (\ref{a6}) and based on the convergence of $T_{l_{k-1}}^0(k-1,N), W_{l_{k-1}}(k-1,N)$, we get that $P_{l_{k-1}}^0(k-1)$ is convergent. Denote $P_{l_1}^0=\lim_{N\rightarrow \infty}P_{l_{k-1}}^0(k-1,N)$, then $P_{l_1}^0$ satisfies the algebraic Riccati equation (\ref{s24}). Now, we show the convergence of $P_{l_{k-1}}(k-1,N)$. For this purpose, we compute the cost $J_N$ with the initial values $u(-1)=u(-2)=\cdots=u(-d)=0$, but
$x_0$ is arbitrary. First, define a Lyapunov function
\begin{eqnarray}
V_N(k,x(k))&=&\mbox{E}\{x(k)'(P_{l_{k-1}}(k-1)-P_{l_{k-1}}^0(k-1))x(k)\nonumber\\
&&-x(k)'\sum_{s=1}^{d-1}(\alpha_{l_{k-1},l_{k-s}}^{d-s}(k-1,k-s))'W_{l_{k-s-1}}(k-s-1)^{-1}\nonumber\\
&&\times \mbox{E}[\alpha_{l_{k-1},l_{k-s}}^{d-s}(k-1,k-s)x(k)|\mathcal{G}_{k-s-1}]\}.\label{a27}
\end{eqnarray}
Applying (\ref{fff1}) and (\ref{a1})-(\ref{a10}), we deduce that
\begin{eqnarray}
&&V_N(k,x(k))-V_N(k+1,x(k+1))\nonumber\\
&=&\mbox{E}\{x(k)'Qx(k)+u(k-d)'Ru(k-d)\nonumber\\
&&-[u(k-d)+W_{l_{k-d}}(k-d)^{-1}\mbox{E}(\alpha_{l_{k-1},l_{k-d+1}}^1(k-1,k-d+1)x(k)|\mathcal{G}_{k-d})]'W_{l_{k-d}}(k-d)\nonumber\\
&&\times [u(k-d)+W_{l_{k-d}}(k-d)^{-1}\mbox{E}(\alpha_{l_{k-1},l_{k-d+1}}^1(k-1,k-d+1)x(k)|\mathcal{G}_{k-d})]\}.\label{a28}
\end{eqnarray}
The expressions (\ref{c1})-(\ref{c4}) play a key role in the derivation of (\ref{a28}). As shown in Remark 5, (\ref{a1})-(\ref{a10}) is also satisfied for $k=d-1,\cdots,0$. So (\ref{a28}) holds for $k=d-1,\cdots,0$. Adding from $k=0$ to $k=N$ on both sides of (\ref{a28}), one yields
\begin{eqnarray}
V_N(0,x_0)&=&\sum_{k=0}^N[V_N(k,x(k))-V_N(k+1,x(k+1))]\nonumber\\
&=&\sum_{k=0}^N\mbox{E}\{x(k)'Qx(k)+u(k-d)'Ru(k-d)\nonumber\\
&&-[u(k-d)+W_{l_{k-d}}(k-d)^{-1}\mbox{E}(\alpha_{l_{k-1},l_{k-d+1}}^1(k-1,k-d+1)x(k)|\mathcal{G}_{k-d})]'\nonumber\\
&&\times W_{l_{k-d}}(k-d)\nonumber\\
&&\times [u(k-d)+W_{l_{k-d}}(k-d)^{-1}\mbox{E}(\alpha_{l_{k-1},l_{k-d+1}}^1(k-1,k-d+1)x(k)|\mathcal{G}_{k-d})]\}.\label{s53}
\end{eqnarray}
It follows from (\ref{s53}) that
\begin{eqnarray}
J_N&=&\mbox{E}\{x_0'(P_{l_{-1}}(-1,N)-P_{l_{-1}}^0(-1,N))x_0+x_0'T_{l_{-1}}^0(-1,N)'
W_{l_{-1}}(-1,N)^{-1}T_{l_{-1}}^0(-1,N)x_0\}\nonumber\\
&=&\mbox{E}\{x_0'P_{l_{-1}}(-1,N)x_0\}.\nonumber
\end{eqnarray}
From the arbitrariness of $x_0$, we get
\begin{eqnarray}
x_0'P_{l_{-1}}(-1,N)x_0=J_N^*\leq J_{N+1}^*=x_0'P_{l_{-1}}(-1,N+1)x_0,\nonumber
\end{eqnarray}
which implies that $P_{l_{-1}}(-1,N)$ increases with respect to $N$. Analogous to the derivation of (73)-(77) in \cite{ZLXF15}, we can show that there exist constant $\lambda$ and $c$, such that
\begin{eqnarray}
J&=&\sum_{k=0}^{\infty}\mbox{E}\{x(k)'Qx(k)\}+\sum_{k=d}^{\infty}\mbox{E}\{u(k-d)'Ru(k-d)\}\nonumber\\
&\leq& 2\lambda cx_0'x_0.\nonumber
\end{eqnarray}
Recalling that $0\leq x_0'P_{l_{-1}}(-1,N)x_0=J_N^*\leq J=2\lambda cx_0'x_0$,\nonumber
which indicates that $0\leq P_{l_{-1}}(-1,N)\leq 2\lambda cI$. The boundedness of $P_{l_{-1}}(-1,N)$ is shown. Recalling that $P_{l_{-1}}(-1,N)$ is monotonically increasing. Therefore, it is convergent. In light of the time invariance of $P_{l_{k-1}}(k-1,N)$, we obtain that
\begin{eqnarray}
\lim_{N\rightarrow \infty}P_{l_{k-1}}(k-1,N+k)=\lim_{N\rightarrow \infty}P_{l_{-1}}(-1,N)=P_{l_{1}}, l_{1}\in \Theta.\nonumber
\end{eqnarray}
This completes the first part of the proof of Theorem 2.

(2) In the second part of the proof, we will show that (\ref{ss29}) is satisfied. Note from Lemma 4, we known that there exists an integer $N_0$, such that
\begin{eqnarray}
&&\Omega_{l_{d-1},l_{d-s}}(d-1,d-s,N_0)\nonumber\\
&=&(P_{l_{d-1}}(d-1,N_0)-P_{l_{d-1}}^0(d-1,N_0))-\sum_{s=1}^{d-1}(\alpha_{l_{d-1},l_{d-s}}^{d-s}(d-1,d-s,N_0))'\nonumber\\
&&\times W_{l_{d-1-s}}(d-1-s,N_0)^{-1}\alpha_{l_{d-1},l_{d-s}}^{d-s}(d-1,d-s,N_0)>0.\nonumber
\end{eqnarray}
On the other hand, $\Omega_{l_{d-1},l_{d-(d-1)}}(d-1,d-(d-1),N)$ is monotonically increasing with respect to $N$, so we have
\begin{eqnarray}
\Omega_{l_{1},l_{d-1}}&=&(P_{l_1}-P_{l_1}^0)-\sum_{s=1}^{d-1}(\alpha_{l_1,l_s}^{d-s})'W_{l_{s+1}}^{-1}\alpha_{l_1,l_s}^{d-s}\nonumber\\
&=&\lim_{N\rightarrow \infty}\Omega_{l_{d-1},l_{d-(d-1)}}(d-1,d-(d-1),N)\nonumber\\
&\geq & \Omega_{l_{d-1},l_{d-(d-1)}}(d-1,d-(d-1),N_0)>0.\nonumber
\end{eqnarray}
Therefore (\ref{ss29}) is satisfied. The proof of Theorem 2 is completed.

\section{Proof of Theorem 3}
\emph{Proof}. Suppose $P_{l_1}, P_{l_1}^0, W_{l_d}$ and $T_{l_d}^j(j=0,1,\cdots,d-1)$ are the solutions to (\ref{s19})-(\ref{s28}) such that
\begin{eqnarray}
(P_{l_1}-P_{l_1}^0)-\sum_{s=1}^{d-1}(\alpha_{l_1,l_s}^{d-s})'W_{l_{s+1}}^{-1}\alpha_{l_1,l_s}^{d-s}>0.\nonumber
\end{eqnarray}
In what follows, we shall prove that the optimal controller (\ref{s57}) stabilizes (\ref{fff1}). For this purpose, we introduce the following Lyapunov function
\begin{eqnarray}
V(k,x(k))&=&\mbox{E}\{x(k)'(P_{l_1}-P_{l_1}^0)x(k)-\sum_{s=1}^{d-1}[x(k)'(\alpha_{l_1,l_s}^{d-s})']W_{l_{s+1}}^{-1}
\mbox{E}[\alpha_{l_1,l_s}^{d-s}x(k)|{\cal{G}}_{k-s-1}]\}.\label{s59}
\end{eqnarray}
Employing (\ref{fff1}) and (\ref{s19})-(\ref{s28}), one gets
\begin{eqnarray}
&&V(k,x(k))-V(k+1,x(k+1))\nonumber\\
&=&\mbox{E}\{x(k)'Qx(k)+u(k-d)'Ru(k-d)\}\geq 0.\label{s60}
\end{eqnarray}
(\ref{s60}) implies that $V(k,x(k))$ decreases with respect to $k$. Further, it follows from (\ref{s59}) that
\begin{eqnarray}
V(k,x(k))
&\geq &
\mbox{E}\{x(k)'[(P_{l_1}-P_{l_1}^0)-\sum_{s=1}^{d-1}(\alpha_{l_1,l_s}^{d-s})'W_{l_{s+1}}^{-1}
\alpha_{l_1,l_s}^{d-s}]x(k)\}\geq 0.\label{s61}
\end{eqnarray}
(\ref{s61}) indicates that $V(k,x(k))$ is bounded below. It follows from (\ref{s60}) and (\ref{s61}) that $V(k,x(k))$ is convergent.

Next, set $m$ to be any nonnegative integer. Via adding from $k=m+d$ to $k=m+N$ on both sides of (\ref{s60}) and letting $m\rightarrow +\infty$, we get that
\begin{eqnarray}
&&\lim_{m\rightarrow\infty} \sum_{k=m+d}^{m+N}\mbox{E}[x(k)'Qx(k)+u(k-d)'Ru(k-d)]\nonumber\\
&=&\lim_{m\rightarrow \infty}V(m+d,x(m+d))-V(m+N+1,x(m+N+1))=0,\label{s62}
\end{eqnarray}
where the last equality holds owning to the convergence of $V(k,x(k))$. Note that
\begin{eqnarray}
&&\sum_{k=d}^N\mbox{E}[x(k)'Qx(k)+u(k-d)'Ru(k-d)]\nonumber\\
&\geq & \mbox{E}\{x(d)'[(P_{l_{d-1}}(d-1,N)-P_{l_{d-1}}^0(d-1,N))-\sum_{s=1}^{d-1}
(\alpha_{l_{d-1},l_{d-s}}^{d-s}(d-1,d-s,N))'\nonumber\\
&&\times (W_{l_{d-s-1}}(d-s-1,N))^{-1}\alpha_{l_{d-1},l_{d-s}}^{d-s}(d-1,d-s,N)]x(d)\}.\nonumber
\end{eqnarray}
By a time shift of length of $m$, it results in
\begin{eqnarray}
&&\sum_{k=m+d}^{m+N}\mbox{E}[x(k)'Qx(k)+u(k-d)'Ru(k-d)]\nonumber\\
&\geq& 
\mbox{E}\{x(m+d)'[(P_{l_{d-1}}(d-1,N)-P_{l_{d-1}}^0(d-1,N))-\sum_{s=1}^{d-1}
(\alpha_{l_{d-1},l_{d-s}}^{d-s}(d-1,d-s,N))'\nonumber\\
&&\times (W_{l_{d-s-1}}(d-s-1,N))^{-1}\alpha_{l_{d-1},l_{d-s}}^{d-s}(d-1,d-s,N)]x(m+d)\}\nonumber\\
&\geq& 0.\label{s63}
\end{eqnarray}
In view of (\ref{s62}), we get that
\begin{eqnarray}
&&\mbox{E}\{x(m+d)'[(P_{l_{d-1}}(d-1,N)-P_{l_{d-1}}^0(d-1,N))-\sum_{s=1}^{d-1}
(\alpha_{l_{d-1},l_{d-s}}^{d-s}(d-1,d-s,N))'\nonumber\\
&&\times (W_{l_{d-s-1}}(d-s-1,N))^{-1}\alpha_{l_{d-1},l_{d-s}}^{d-s}(d-1,d-s,N)]x(m+d)\}\nonumber\\
&=& 0.\label{s64}
\end{eqnarray}
Recalling from Lemma 4, we know that there exists an integer $N_0$, such that
\begin{eqnarray}
&&(P_{l_{d-1}}(d-1,N_0)-P_{l_{d-1}}^0(d-1,N_0))-\sum_{s=1}^{d-1}
(\alpha_{l_{d-1},l_{d-s}}^{d-s}(d-1,d-s,N_0))'\nonumber\\
&&\times (W_{l_{d-s-1}}(d-s-1,N_0))^{-1}\alpha_{l_{d-1},l_{d-s}}^{d-s}(d-1,d-s,N_0)>0.\nonumber
\end{eqnarray}
So (\ref{s64}) indicates that $\lim_{m\rightarrow \infty}\mbox{E}[x(m+d)'x(m+d)]=0$. That is to say (\ref{s57}) stabilizes (\ref{fff1}) in the mean-square sense.

In what follows, we will show that the cost function (\ref{ss1}) is minimized by (\ref{s57}). Adding from $k=0$ to $k=N$ to (\ref{s60}), one yields
\begin{eqnarray}
&&\mbox{E}\{\sum_{k=0}^Nx(k)'Qx(k)+\sum_{k=d}^Nu(k-d)'Ru(k-d)\}\nonumber\\
&=&V(0,x_0)-V(N+1,x(N+1))\nonumber\\
&&+\sum_{k=d}^N\mbox{E}\{[u(k-d)+W_{l_d}^{-1}T_{l_d}^0x(k-d+1)+\sum_{j=1}^{d-1}W_{l_d}^{-1}T_{l_d}^ju(k-2d+j)]'\nonumber\\
&&\times W_{l_d}[u(k-d)+W_{l_d}^{-1}T_{l_d}^0x(k-d+1)+\sum_{j=1}^{d-1}W_{l_d}^{-1}T_{l_d}^ju(k-2d+j)]\}\nonumber\\
&&+\sum_{k=0}^{d-1}\mbox{E}\{[u(k-d)+W_{l_d}^{-1}T_{l_d}^0x(k-d+1)+\sum_{j=1}^{d-1}W_{l_d}^{-1}T_{l_d}^ju(k-2d+j)]'\nonumber\\
&&\times W_{l_d}[u(k-d)+W_{l_d}^{-1}T_{l_d}^0x(k-d+1)+\sum_{j=1}^{d-1}W_{l_d}^{-1}T_{l_d}^ju(k-2d+j)]\}\nonumber\\
&&-\sum_{k=0}^{d-1}\mbox{E}\{u(k-d)'Ru(k-d)\}.\label{s65}
\end{eqnarray}
Note that system (\ref{fff1}) is mean square stabilizable, we have
\begin{eqnarray}
\lim_{k\rightarrow \infty}\mbox{E}\{x(k)'P_{l_1}x(k)\}=0.\label{sss65}
\end{eqnarray}
On the other hand,
\begin{eqnarray}
0&\leq &V(k,x(k))\nonumber\\
&=&\mbox{E}\{x(k)'(P_{l_1}-P_{l_1}^0)x(k)-\sum_{s=1}^{d-1}[x(k)'(\alpha_{l_1,l_s}^{d-s})']
W_{l_{s+1}}^{-1}\mbox{E}[\alpha_{l_1,l_s}^{d-s}x(k)|{\cal{G}}_{k-s-1}]\}\nonumber\\
&=&\mbox{E}\{x(k)'(P_{l_1}-P_{l_1}^0)x(k)-\sum_{s=1}^{d-1}\mbox{E}[x(k)'(\alpha_{l_1,l_s}^{d-s})'|{\cal{G}}_{k-s-1}]
W_{l_{s+1}}^{-1}\mbox{E}[\alpha_{l_1,l_s}^{d-s}x(k)|{\cal{G}}_{k-s-1}]\}\nonumber\\
&\leq &\mbox{E}\{x(k)'P_{l_1}x(k)\}.\label{ss65}
\end{eqnarray}
It concludes from (\ref{sss65}) and (\ref{ss65}) that
\begin{eqnarray}
\lim_{k\rightarrow \infty}V(k,x(k))=0.\nonumber
\end{eqnarray}
Let $N\rightarrow \infty$ on both sides of (\ref{s65}), we get
\begin{eqnarray}
J&=&V(0,x_0)-\sum_{k=0}^{d-1}\mbox{E}\{u(k-d)'Ru(k-d)\}\nonumber\\
&&+\sum_{k=d}^N\mbox{E}\{[u(k-d)+W_{l_d}^{-1}T_{l_d}^0x(k-d+1)+\sum_{j=1}^{d-1}W_{l_d}^{-1}T_{l_d}^ju(k-2d+j)]'\nonumber\\
&&\times W_{l_d}[u(k-d)+W_{l_d}^{-1}T_{l_d}^0x(k-d+1)+\sum_{j=1}^{d-1}W_{l_d}^{-1}T_{l_d}^ju(k-2d+j)]\}\nonumber\\
&&+\sum_{k=0}^{d-1}\mbox{E}\{[u(k-d)+W_{l_d}^{-1}T_{l_d}^0x(k-d+1)+\sum_{j=1}^{d-1}W_{l_d}^{-1}T_{l_d}^ju(k-2d+j)]'\nonumber\\
&&\times W_{l_d}[u(k-d)+W_{l_d}^{-1}T_{l_d}^0x(k-d+1)+\sum_{j=1}^{d-1}W_{l_d}^{-1}T_{l_d}^ju(k-2d+j)]\}.\label{s66}
\end{eqnarray}
Since $W_{l_d}$ is positive, then (\ref{s66}) is minimized if and only if
\begin{eqnarray}
u(k-d)=-W_{l_d}^{-1}T_{l_d}^0x(k-d+1)-\sum_{j=1}^{d-1}W_{l_d}^{-1}T_{l_d}^ju(k-2d+j).\nonumber
\end{eqnarray}
The corresponding optimal cost (\ref{s58}) is obtained.

From the derivation of Theorem 2, the following condition is satisfied
\begin{eqnarray}
(P_{l_1}-P_{l_1}^0)-\sum_{s=1}^{d-1}(\alpha_{l_1,l_s}^{d-s})'W_{l_{s+1}}^{-1}\alpha_{l_1,l_s}^{d-s}>0.\nonumber
\end{eqnarray}
The uniqueness can be similarly shown as the proof of Theorem 3 in \cite{ZLXF15}, so we omit here. This completes the proof of Theorem 3.


\end{document}